# A7: an aperiodic set of 7 square dominoes

Vincent Van Dongen

Independent researcher, Canada; vincent.vandongen@gmail.com

**Abstract**

This paper presents an aperiodic tileset of 7 square dominoes. We call it A7 as it directly relates to the aperiodic set Ammann A3. We start with a description of the tileset. We then present Ammann A3 and its direct link with tileset A7.

## Introduction

Wang tiles are square dominoes that can't rotate or reflect [1]. Ten years ago, an aperiodic set of 11 tiles was published [2], proven to be the one of minimal size. In this document, we explore a new set of square dominoes. These are not Wang tiles as they can rotate. Other work on square tiles with marks for local constraints can be found in the literature such as [4][5]. In this paper, we present an aperiodic set of 7 square dominoes. The set is directly based on the aperiodic set Ammann A3 [6]. We first present A7 and then present its link with A3.

## A tileset of 7 square dominoes

Tileset A7 is a new set of tiles that consists of 7 square dominoes. Each border is identified by a pair of colors. A tiler is allowed to rotate the tiles to position them in the desired orientation but border colors of adjacent tiles must match.

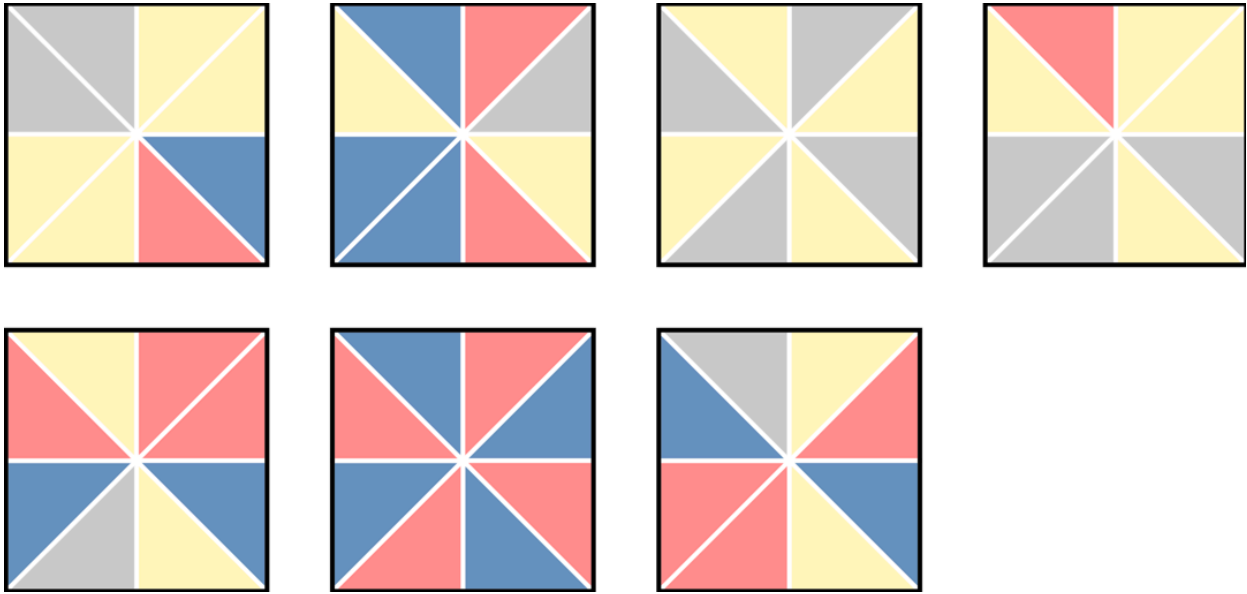

***Figure 1:*** *The 7 square dominoes of A7 with each border identified with two colors.*

Figure 2 gives an example of an A7 tiling of size 14x22, i.e. 308 square dominoes.

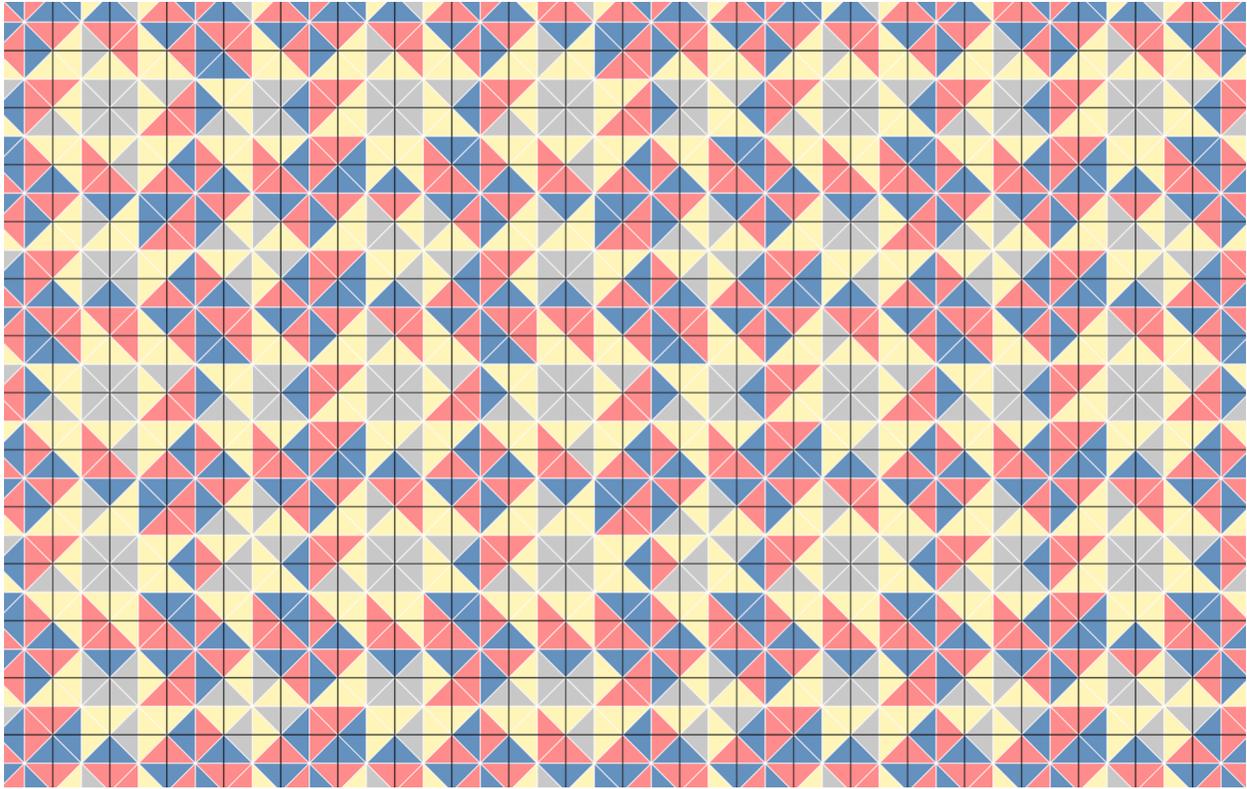

*Figure 2: Example of a 14x22 tiling A7.*

## Aperiodic set Ammann A3

The aperiodic set Ammann A3 consists of 3 tiles. See Figure 3. Marks on the tiles enforce non-periodic tilings. For a tiling to be valid, lines must be continuous. These marks can be implemented in a physical manner with borders with etches where the white lines are. See Figure 4.

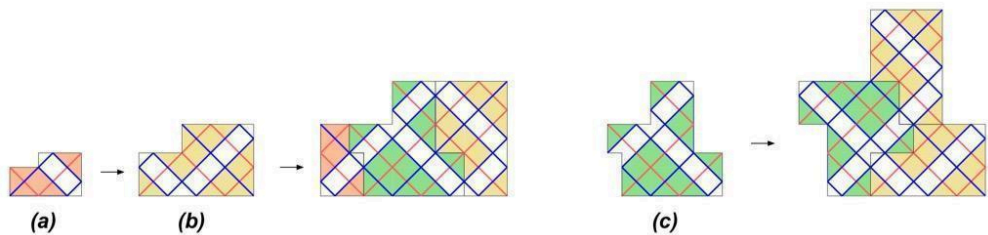

*Figure 3: Aperiodic set Ammann A3 consists of 3 tiles with marks that enforce non-periodic tiling.*

These marks can be implemented in a physical manner with borders with etches where the white lines are. See Figure 4.

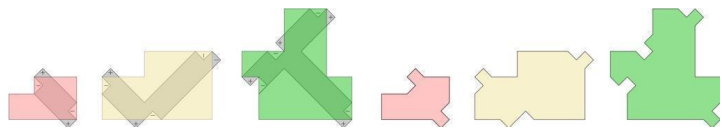

*Figure 4: Aperiodic set Ammann A3 implemented with borders with etches.*

We carefully analysed these marks on large tilings and we noticed that the entire tiling can be created with 7 supertiles of rectangular shapes. Figure 5 shows these supertiles within A3 tilings.

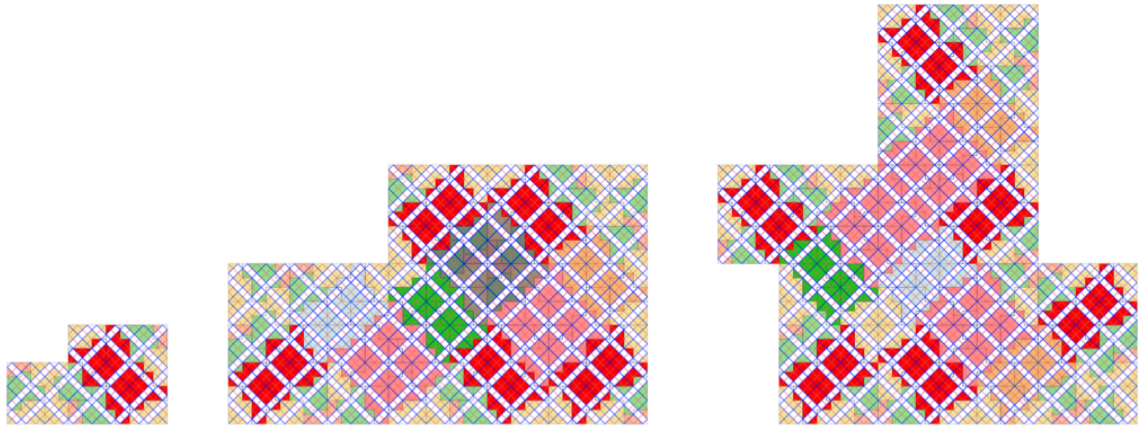

*Figure 5: Creation of supertiles that are of rectangular shape.*

This gave us the idea of associating a square domino to each of them. Figure 6 shows the 7 dominoes associated with the supertiles. We associated a pair of colors for every border shape. This way, only when border colors match on the dominoes, the supertiles can be assembled with no overlap and no hole. A key observation is that the set of 7 supertiles is aperiodic based on the fact that their shapes are such that when border to border properly meet, the white lines are automatically respected. And it was already proven that set A3 is aperiodic.

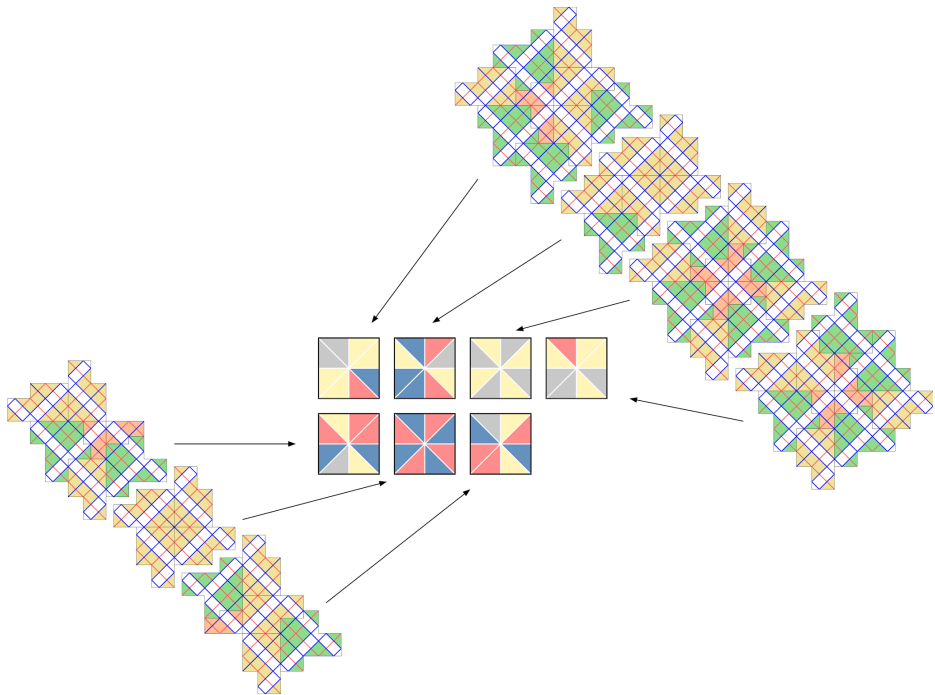

*Figure 6: The 7 supertiles of rectangular shapes are mapped to square dominoes.*

## Rules for tiling with A7

How can we tile with A7? One way is by trials and errors but this is quite challenging. A more efficient approach is to make use of iterative rules. Fortunately we could adapt the ones of A3. The rules for A7 are shown in Figure 7. The starting points are 3 supertiles of A7 square dominoes. See Figure 8.

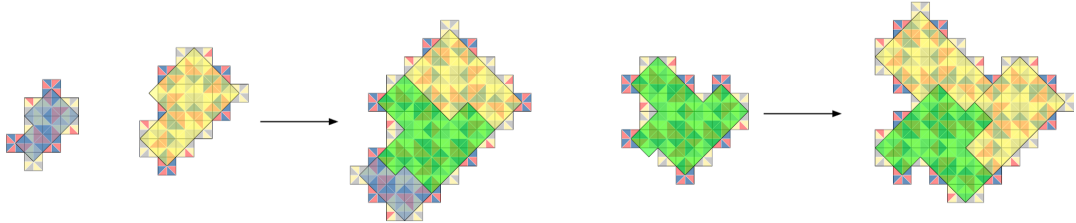

*Figure 7: Just like for A3, two rules can be used for creating A7 tilings.*

The shapes on Figure 8 contain overlaps. We carefully analysed these overlaps to make sure that they are always duplicate A7 dominoes that can be removed.

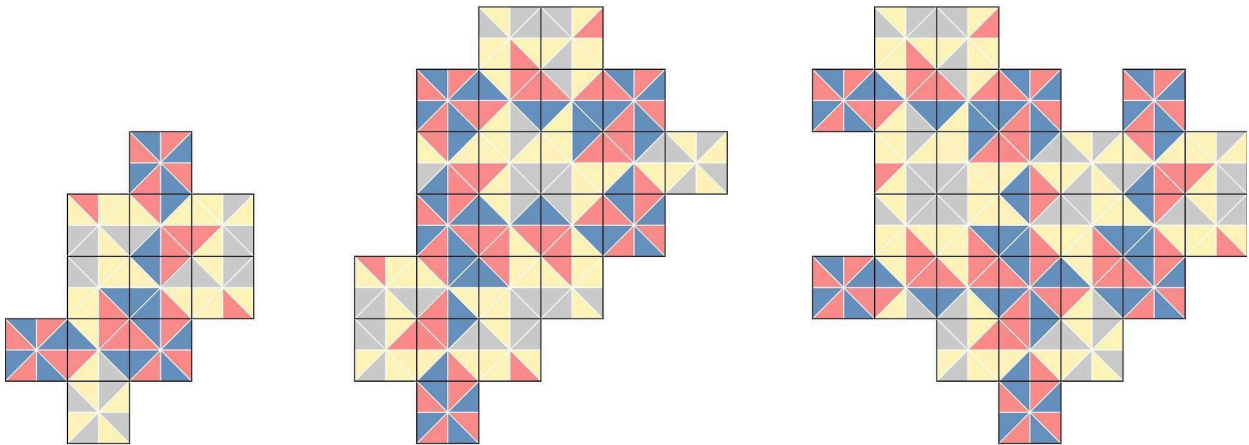

*Figure 8: Three tilings made of A7 dominoes that serve as rules for creating larger tilings.*

## Fixed points for A7

The 3 rules provided above allow us to create tilings of any size. The question now is: Can we have tiling rules with no overlap? The topic of this section is to provide such rules. In this approach, a rule is associated with each domino. This means 7 rules, one rule per square domino. What we know is that A7 tiling can be marked with a grid of Ammann. This grid subdivides A7 tiling into supertiles that correspond to the fixed points of the dominoes. Figure 9 shows the set of non-overlapping supertiles associated to A7. The square dominoes are here labelled 0, 1, 2, 3, A, B, C. This helps to visualize the use of the square dominoes in the supertiles. The corners of the supertiles are always one of the two tiles A or B. The Ammann grid lies at the boundary of the supertiles. The reason for calling these supertiles fixed points is as follows: the square domino tile borders keep their 'color type' at the supertile level. In other words, if a square domino has a border that is compatible with other tile borders, at the supertile level, its corresponding border offers the same border-to-border compatibility. And this will be true at every iteration.

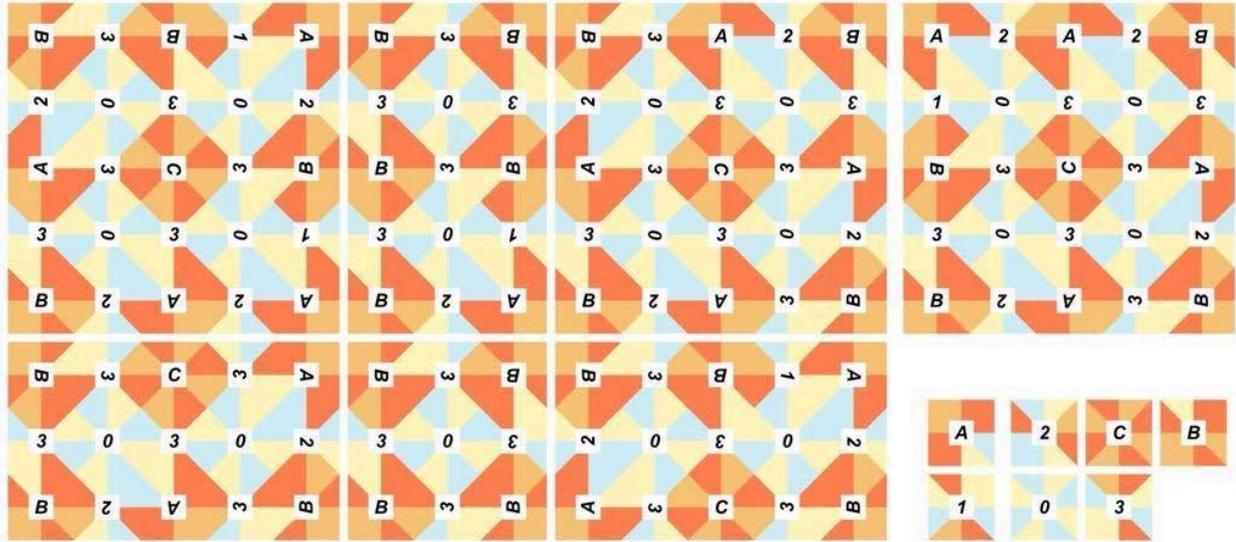

*Figure 9:* *Supertiles of A7 square dominoes.*

## Conclusion

Tileset A7 is a set of 7 square dominoes based on Ammann A3. This set was obtained by creating rectangular supertiles of A3 tiles. We showed that rules based on A3 can be used to tile with these dominoes. A similar approach than the one briefly presented in this paper could be taken and applied to other aperiodic tilings. Now that we found this set, we conclude with an open question: is there an aperiodic set of square dominoes of size smaller than 7?

## Acknowledgement

Special thanks to Pierre Gradit and Jean-François Husson for their precious comments and great support.

## References

[1] Wang tile on Wikipedia. https://en.wikipedia.org/wiki/Wang_tile
[2] Jeandel, E. and Rao, M., "An aperiodic set of 11 Wang tiles.". 2015. https://arxiv.org/abs/1506.06492
[3] Van Dongen, V. and Gradit, P., 2024, https://arxiv.org/pdf/2407.06202
[4] Yliopisto T., "Robinson's aperiodic tileset", 2021, https://users.utu.fi/jkari/wp-content/uploads/sites/1251/2021/10/robinson.pdf
[5] Lagae, A.; Kari, J. and Dutré, P., "Aperiodic Sets of Square Tiles with Colored Corners.", 2006, KUL Report CW460.
[6] Ammann, R., Grunbaum B. and Shephard G. C. "Aperiodic tiles." Discrete & Computational Geometry, vol. 8, 1992.